\documentclass[reqno,11pt]{amsart}
\usepackage{amsmath,amssymb,graphics,epsfig,color,cite}
\newtheorem{theorem}{Theorem}[section]
\newtheorem{proposition}[theorem]{Proposition}

\theoremstyle{remark}
\newtheorem{remark}[theorem]{Remark}
\theoremstyle{definition}

\numberwithin{equation}{section}
\numberwithin{theorem}{section}

\newcommand{\mc}[1]{{\mathcal #1}}
\newcommand{\bb}[1]{{\mathbb #1}}

\newcommand{\rme}{\mathrm{e}}
\newcommand{\rmi}{\mathrm{i}\,}

\newcommand{\conj}[1]{\overline{#1}}

\newcommand{\lan}{\left\langle}
\newcommand{\ran}{\right\rangle}

\begin{document}

\title{Structured conditioning of Hamiltonian eigenvalue problems}

\author [P.\ Butt\`a] {Paolo Butt\`a}
\address{Paolo Butt\`a, Dipartimento di Matematica, SAPIENZA Universit\`a di Roma, P.le Aldo Moro 5, 00185 Roma, Italy}
\email{butta@mat.uniroma1.it}

\author [S.\ Noschese] {Silvia Noschese}
\address{Silvia Noschese, Dipartimento di Matematica, SAPIENZA Universit\`a di Ro\-ma, P.le Aldo Moro 5, 00185 Roma, Italy}
\email{noschese@mat.uniroma1.it}



\begin{abstract}
We discuss the effect of structure-preserving perturbations on complex or real Hamiltonian eigenproblems and characterize the structured worst-case effect perturbations.  We derive significant expressions  for both the structured condition numbers and  the worst-case effect Hamiltonian perturbations. It is shown that, for purely imaginary eigenvalues, the usual unstructured perturbation analysis is sufficient.
\end{abstract}

\maketitle
\thispagestyle{empty}

\section{Introduction}
\label{sec:1}

For structured eigenvalue problems, algorithms that preserve the underlying matrix structure may improve the accuracy and efficiency of the eigenvalue computation and preserve possible eigenvalue symmetries in finite-precision arithmetic. The concept of  backward stability, i.e., the requirement that the computed eigenvalues are the exact eigenvalues of a slightly perturbed matrix, is extended by requiring the perturbed matrix  to have the same structure as the original one \cite{B87}. In order to assess the {\it strong} numerical backward-stability of the employed structured algorithm and not to overestimate the worst-case effect of  perturbations, it is appropriate to consider suitable measures of the sensitivity of the eigenvalues to perturbations  of the same structure. 

Our analysis is based on the perturbation expansion of a simple eigenvalue. The structured condition number of an eigenvalue $\lambda \in A$ is indeed a first-order measure of the worst-case effect  on $\lambda$ of perturbations of the same structure as  $A$. The structured conditioning measures we deal with can be computed endowing the subspace of matrices  with the Frobenius norm; see, e.g., \cite{T03, KKT, R06} and references therein.
 
Here we are concerned with the Hamiltonian eigenvalue problems. They arise from a number of applications, particularly in systems and control theory. Algorithms and applications of Hamiltonian eigenproblems are discussed, e.g., in \cite{BKM05} and references therein.

We investigate  the sensitivity of the eigenvalues of a complex [real] Hamiltonian matrix with respect to  complex [real] Hamiltonian perturbations.   

Attention has been paid  to derive straightforwardly computable formulae, for both  the structured condition numbers and  the maximal structured perturbations, that result to be explicative of the relationship between unstructured and structured conditioning. 
It is well known,  see, e.g.,  \cite{W65}, that the worst (unstructured) perturbation which may affect a simple eigenvalue $\lambda $ of a given matrix $A$ arises under the action of the matrix $yx^{*}$, with $x$ and $y$ corresponding right and left eigenvectors, normalized to have $\left\| x\right\| _{2}=\left\|y\right\| _{2}=1$.  As a matter of  fact, if we add the Hamiltonian-structure requirement,  just as for sparsity-structures \cite{NP06} or symmetry-patterns \cite{NP07}, we prove that  the Hamiltonian matrix which yields the structured worst-case effect perturbation is the {\it structured analogue} of  $yx^{*}$, i.e. the normalized projection of  $yx^{*}$ onto the space of the complex [real] Hamiltonian matrices. Notice that  a further condition on the angle between the right and left eigenvectors is here required in the normalization.

The paper is organized as follows. Section \ref{sec:2} provides formulae for the eigenvalue structured condition numbers and  for the maximal Hamiltonian perturbations. In Section \ref{sec:3}, the real case is addressed. Numerical tests are presented in Section \ref{sec:4}. Finally, conclusions are drawn in Section \ref{sec:5}.

\section{Eigenvalue structured conditioning of Hamiltonian matrices}
\label{sec:2}

We start with some notation and definitions. We denote by ${\mc H}_{\bb C}$ the linear subspace of $2n$-dimensional complex Hamiltonian matrices, i.e., 
\[
\begin{split}
\mc H_{\bb C} & = \left\{Q \in \bb C^{2n\times 2n} \colon QJ=(QJ)^*\right\} \\ & = \left\{Q=\begin{pmatrix} K & M \\ L & -K^* \end{pmatrix} \colon K,L,M \in \bb C^{n\times n}\;,\;L=L^*\;,\;M=M^*\right\}\;,
\end{split}
\]
where $J$ is the fundamental symplectic matrix, i.e., 
\[
J=\begin{pmatrix} 0 & I_n \\ -I_n & 0 \end{pmatrix}\;,
\]
with $I_n$ the $n \times n$ identity matrix. We observe that the Hamiltonian matrix $Q$ is uniquely determined by the $4n^2$ real parameters given by the $2n^2$ real and imaginary parts of the elements of $K$, the $2n$ (real) diagonal elements of $L,M$, and the $2n(n-1)$ real and imaginary parts of the elements of the strict lower triangular parts of $L,M$.
\begin{proposition}
\label{prop:1}
The closest Hamiltonian matrix to a given matrix $A = \begin{pmatrix} A_1 & A_3 \\ A_2 & A_4 \end{pmatrix}\in \bb C^{2n\times 2n}$ with respect to the Frobenius norm is 
\[
A|_{\mc H_{\bb C}} = \frac 12 \begin{pmatrix} A_1-A_4^* & A_3+A_3^* \\ A_2+A_2^* & A_4-A_1^* \end{pmatrix} = \frac 12 (A+JA^* J)\;.
\]
\end{proposition}
\proof For $Q = \begin{pmatrix} K & M \\ L & -K^* \end{pmatrix} \in \mc H$ we have
\[
\|A-Q\|_F^2=\|K-A_1\|_F^2 + \|K+A_4^*\|_F^2 +\|L-A_2\|_F^2+\|M-A_3\|_F^2\;.
\]
Since $L,M$ are only requested to be Hermitian, it is known that the last two terms in the right-hand side are minimized for $L=\frac 12(A_2+A_2^*)$ and $M = \frac 12 (A_3+A_3^*)$. On the other hand, the minimization of
\[
G(K) = \|K-A_1\|_F^2 + \|K+A_4^*\|_F^2
\]
leads to $K=\frac 12 (A_1-A_4^*)$. The proposition is thus proved.
\qed

\medskip
To state the next results we need some additional notation and definitions. It will be useful the following characterization of the space of complex Hamiltonian matrices,
\[
{\mc H}_{\bb C}  = \left\{H+\rmi W \colon H \in \mc H, \;\; W\in \mc W\right\}\;,
\]
where $\mc H$ [resp.\ $\mc W]$ is the linear space of real Hamiltonian [resp.\ skew-Hamiltonian] matrices, i.e.,
\[
\begin{split}
\mc H & = \left\{H \in \bb R^{2n\times 2n} \colon HJ=(HJ)^T\right\}\;, \\ \mc W & = \left\{W \in \bb R^{2n\times 2n} \colon WJ=-(WJ)^T\right\}\;.
\end{split}
\]
In particular, if $A=\Re(A)+\rm i \Im(A)$ then $A|_{\mc H_{\bb C}} = \Re(A)|_{\mc H} +\rm i \Im(A)|_{\mc W}$ where,  for $B\in \bb R^{2n\times 2n}$,
\[
B|_{\mc H} := \frac 12 (B+JB^TJ)\;\,\qquad B|_{\mc W} := \frac 12 (B-JB^TJ)\;.
\]
We finally introduce the normalized projection of a matrix $A$ onto $\mc H_{\bb C}$ as
\begin{equation}
\label{p:0}
A|_\mc N := \frac{A|_{\mc H_{\bb C}}}{\|A|_{\mc H_{\bb C}}\|_F}\;.
\end{equation}
\begin{theorem}
\label{teo:1}
Let $\lambda$ be a simple eigenvalue of a Hamiltonian matrix $Q$, with corresponding right and left eigenvectors $x$ and $y$ normalized to have
\begin{equation}
\label{p:1}
\|x\|_2=\|y\|_2 = 1\;, \qquad \Im(y^*Jx) = 0\;.
\end{equation}
Given any Hamiltonian matrix $E$ with $\|E\|_F = 1$, let $\lambda_E(t)$ be an eigenvalue of $Q+tE$ converging to $\lambda$ as $t\to 0$. Then,
\begin{equation}
\label{p:2-}
|\dot\lambda_E(0)| \le \max\left\{\left|\frac{y^*Gx}{y^*x}\right| \colon  \|G\|_F = 1,\, G \in \mc H_{\bb C}\,\right\} = \frac{\|yx^*|_{\mc H_{\bb C}}\|_F}{|y^*x|}\;.
\end{equation}
Moreover, 
\begin{equation}
\label{p:2}
\dot\lambda_E(0) = \frac{\|yx^*|_{\mc H_{\bb C}}\|_F}{y^*x} \qquad \mathrm{if} \qquad E=\pm yx^*|_{\mc N}\;.
\end{equation}
\end{theorem}
\begin{remark}
\label{rem:p1}
By standard perturbation theory, the variational estimate in \eqref{p:2-} holds for any choice of the right and left eigenvectors $x$ and $y$. Instead, the assumption \eqref{p:1} is needed for the explicit expressions of the maximum and of the maximizers to hold. 
\end{remark}
\noindent
{\it Proof of Theorem \ref{teo:1}}. 
We shall use the real scalar product $\lan A,B\ran := \mathrm{Trace}(A^TB)$, $A,B\in \bb R^{2n\times 2n}$, note $\|A\|_F=\sqrt{\lan A,A\ran}$. Our target is to compute
\[
\mathrm{argmax}\left\{|y^*(H+\rmi W) x| \colon \;H\in\mc H\;,\; W\in\mc W\;,\;\|H\|_F^2+\|W\|_F^2=1\right\}.
\]
The requirement $H\in \mc H$ means that $HJ$ is symmetric, i.e., it is orthogonal (w.r.t.\ the above scalar product) to the space of skew-symmetric matrices. Similarly, $W\in \mc W$ means that $WJ$ is orthogonal to the space of symmetric matrices. Therefore, we study the auxiliary (unconstrained)  variational problem for the function
\[
G(H,W) = |y^*(H+\rmi W) x| ^2+ \lambda (\|H\|_F^2+\|W\|_F^2-1) + \lan\Lambda,HJ\ran + \lan\Gamma,WJ\ran\;,
\]
which depends on $4n^2+1$ Lagrange multipliers, precisely the real parameter $\lambda$, the skew-symmetric matrix $\Lambda$, and  the symmetric matrix $\Gamma$. Such parameters will be eventually determined by requiring that the solution satisfies all the constraints. 

It is now useful to introduce the matrices $\xi = \Re(yx^*)$,  $\eta = \Im(yx^*)$, so that
\[
yx^* = \xi + \rmi \eta\;,\quad \bar y x^T  = \xi - \rmi \eta
\]
and therefore
\[
\begin{split}
 y^*(H+\rmi W) x & = \mathrm{Trace}((\bar y x^T)^T(H+\rmi W)) = \mathrm{Trace}((\xi^T-\rmi\eta^T ) (H+\rmi W)) \\ & = \lan \xi,H\ran +\lan\eta,W\ran + \rmi (\lan \xi,W\ran - \lan \eta, H\ran)\;.
 \end{split}
\]
Whence,
\[
\begin{split}
G(H,W) & = (\lan \xi,H\ran +\lan\eta,W\ran)^2+ (\lan \xi,W\ran - \lan \eta, H\ran)^2 \\ & \quad + \lambda (\|H\|_F^2+\|W\|_F^2-1) + \lan\Lambda,HJ\ran + \lan\Gamma,WJ\ran\;,
\end{split}
\]
and, by \eqref{p:1}, also noticing  that $y^*Jx = \mathrm{Trace}(xy^*J) = \lan\xi,J\ran -\rmi \lan\eta,J\ran$, 
\begin{equation}
\label{p:1bis}
\|\xi\|_F^2+\|\eta\|_F^2 = 1\;, \qquad \lan \eta,J\ran = 0\;.
\end{equation}
We now look for the critical points of $G$. The gradients $\nabla_HG(H,W)$, $\nabla_WG(H,W)$, expressed in the form of $2n\times 2n$ matrices, read,
\[
\begin{split}
\nabla_HG(H,W) & = 2a \xi- 2b \eta +2\lambda H + \Lambda J^T\;, \\ \nabla_WG(H,W) & = 2a\eta +2b\xi+2\lambda W + \Gamma J^T\;,
\end{split}
\]
where we introduced the scalar functions,
\[
a = \lan \xi,H\ran +\lan\eta,W\ran\;, \qquad  b=\lan \xi,W\ran - \lan \eta, H\ran\;.
\]
The critical points $(H,W)$ are therefore solutions to
\begin{equation}
\label{HW}
\begin{split}
H & = - \frac a\lambda \xi + \frac b \lambda \eta  - \frac 1{2\lambda} \Lambda J^T\;, \\ W & = - \frac a\lambda \eta  - \frac b\lambda \xi - \frac 1{2\lambda} \Gamma J^T\;.
\end{split}
\end{equation}
By imposing $HJ=(HJ)^T$ and $WJ=-(WJ)^T$, using that $\Lambda=-\Lambda^T$, $\Gamma=\Gamma^T$, and recalling that $JJ^T=I$, $J^T=-J$, we get,
\[
\Lambda = - a (\xi J+J\xi^T) + b (\eta J +J\eta^T) \;,\quad
\Gamma = - a (\eta J-J\eta^T) - b (\xi J -J\xi^T)\;.
\]
Plugging these expressions in \eqref{HW} and calling $\alpha=-\lambda^{-1}a$, $\beta= \lambda^{-1}b$, we conclude that the critical points have to be of the following form,
\begin{equation}
\label{ab}
\begin{split}
H & = \alpha (\xi+J\xi^TJ)  + \beta (\eta + J\eta^TJ) \;, \\ W & = \alpha (\eta-J\eta^TJ)  - \beta (\xi - J\xi^TJ) \;.
\end{split}
\end{equation}
To determine the (still) unknown real parameters $(\alpha,\beta)$, it remains to maximize $|y^*(H+\rmi W) x|^2$ with the constraint $\|H\|_F^2+\|W\|_F^2=1$ for $(H,W)$ as in \eqref{ab}. To compute these quantities we first notice that, as $(\xi^T \pm \rmi \eta^T) J$ are rank-one (complex) matrices, 
\[
\mathrm{Trace}((\xi^T+\rmi\eta^T) J (\xi^T\pm\rmi\eta^T)J)  = \mathrm{Trace}((\xi^T+\rmi\eta^T) J) \mathrm{Trace}(\xi^T\pm\rmi\eta^T)J)\;,
\]
from which, by identifying the real and imaginary parts, we get
\[ 
\begin{split}
\mathrm{Trace}(\xi^TJ\xi^TJ) & = \mathrm{Trace}(\xi^TJ)^2\;,\quad  \mathrm{Trace}(\eta^TJ\eta^TJ) = \mathrm{Trace}(\eta^TJ)^2\;, \\ \mathrm{Trace}(\xi^TJ\eta^TJ) & = \mathrm{Trace}(\eta^TJ\xi^TJ)  = \mathrm{Trace}(\eta^TJ) \mathrm{Trace}(\xi^TJ)\;.
\end{split}
\]
Therefore, by \eqref{p:1bis},
\[
\begin{split}
\lan\xi,J\xi^TJ\ran & = \lan \xi , J\ran^2\;, \quad \lan\eta,J\eta^TJ\ran  = \lan \eta ,J\ran^2=0\;, \\ \lan\xi,J\eta^TJ\ran & = \lan\eta,J\xi^TJ\ran = \lan \xi,J \ran \lan \eta,J\ran=0\;,
\end{split}
\]
so that, one easily computes,
\[
| y^*(H+\rmi W) x|^2 =(1+\lan\xi,J\ran^2)^2  \alpha^2 + (1-\lan\xi,J\ran^2)^2\beta^2\;,
\]
\begin{equation}
\label{nh}
\|H\|_F^2+\|W\|_F^2= 2(1+\lan\xi,J\ran^2)\alpha^2 + 2(1-\lan\xi,J\ran^2)\beta^2\;.
\end{equation}
Moreover, as $yx^*|_{\mc H}$ [resp.\ $(\rmi yx)^*|_{\mc H}$] is given by \eqref{ab} for $(\alpha,\beta) = (1/2,0)$ [resp.\ $(\alpha,\beta)=(0,1/2)$], by \eqref{nh} we get
\[
\|yx^*|_{\mc H}\|_F^2 = \frac{1+\lan\xi,J\ran^2}2\;, \quad \|(\rmi yx)^*|_{\mc H}\|_F^2 = \frac{1-\lan\xi,J\ran^2}2\;.
\]
If $\|(\rmi yx)^*|_{\mc H}\| >0$, i.e., $|\lan\xi,J\ran|< 1$, the (constrained) maximum of $| y^*(H+\rmi W) x|^2$ is reached for $\alpha^2= [2(1+\lan\xi,J\ran^2)]^{-1}$ and $\beta^2 =0$, so that the maximal perturbation is obtained for 
\[
E_\pm = \pm \frac{ \xi+J\xi^TJ}{\sqrt{2(1+\lan\xi,J\ran^2)}} \pm\rmi \frac{ \eta-J\eta^TJ}{\sqrt{2(1+\lan\xi,J\ran^2)}} = \pm yx^*|_{\mc N}\;.
\]
On the other hand, if $\|(\rmi yx)^*|_{\mc H}\| =0$ then $\alpha^2=[2(1+\lan\xi,J\ran^2)]^{-1}$ and $H+\rmi W = E_\pm $ holds trivially. In both cases,

\[
| y^*E_\pm x| = \sqrt{\frac{1+\lan\xi,J\ran^2}2} = \|yx^*|_{\mc H}\|_F\;.
\]
The theorem is thus proved.
\qed
\begin{remark}
\label{rem:p2}
By Theorem \ref{teo:1} the structured condition number is explicitly given by
\[
\kappa_{\mc H_{\bb C}}(\lambda) = \frac{\|yx^*|_{\mc H}\|_F}{|y^*x|}\;,
\]
which can be also used for the structured backward error analysis of the algorithms in \cite{BMX98, BK05}. Since the (unstructured) condition number of a simple eigenvalue is $\kappa(\lambda) = \|yx^*\|_F/|y^*x|$, we have $\kappa_{\mc H_{\bb C}}(\lambda) = \kappa(\lambda)$ if $yx^*$ is Hamiltonian. As shown in the next proposition, this can occur only if $\lambda$ is a purely imaginary eigenvalue.
\end{remark}
\begin{proposition}
\label{prop:p1}
Let $\lambda$ be a simple eigenvalue of a Hamiltonian matrix $Q$ with corresponding right and left eigenvectors $x$ and $y$ normalized such that $\|x\|_2=\|y\|_2 = 1$. If $yx^*$ is a Hamiltonian matrix or a skew-Hamiltonian matrix (i.e., $(yx^*J)=-(yx^*J)^*$) then $\Re\lambda=0$.
\end{proposition}
\proof We prove the statement in the case $yx^*$ is a Hamiltonian matrix, the case $yx^*$ is  skew-Hamiltonian can be treated similarly. We first show that if $yx^*$ is Hamiltonian then $y=\pm Jx$. Indeed, since $J^T=-J$, the identity $yx^*J=(yx^*J)^*$ can be recast as $y_k (J\bar x)_j = (Jx)_k\bar y_j$ for any $k,j=1,\ldots,2n$. As $y^*x\ne 0$, there exists at least one index $k$ for which $(Jx)_k, y_k \ne 0$. Then, it is readily seen that the above identity implies $Jx =c y$ for some $c\in \bb R$. As $\|x\|_2=\|y\|_2 = 1$, we conclude that $y=\pm Jx$. It remains to show that $y=\pm Jx$ implies $\Re\lambda =0$. To this end, it is enough to notice that, as $Q$ is Hamiltonian and $J^2=-1$, the identity $Qx=\lambda x$ can be recast as $(Jx)^*Q= -\bar\lambda (Jx)^*$. Since also $y=\pm Jx$ is solution to $y^*Q=\lambda y^*$, we conclude that $\lambda+\bar\lambda=0$.  
\qed

\section{Structured conditioning of real Hamiltonian eigenproblems}
\label{sec:3}

In several situations, we are interested in real structured perturbations of real Hamiltonian matrices. For these purposes, we are led to consider the real structured condition number $\kappa_{\mc H}(\lambda)$, where only real Hamiltonian perturbations are taken into account. For any smooth manifold $\mc S$, an explicit formula for the real structured condition number $\kappa_{\mc S}(\lambda)$ is given in \cite[Eq.\ (2.16)]{KKT}. Unfortunately, such formula, as its counterpart in the complex case, requires the characterization of the tangent space of $\mc S$ and the construction of a pattern matrix relevant to $\mc S$. On the other hand, it yields the following bound,
\[
\frac{\kappa_{\mc H_{\bb C}}(\lambda)}{\sqrt 2} \le \kappa_{\mc H}(\lambda)\le  \kappa_{\mc H_{\bb C}}(\lambda)\;. 
\]
Actually, we can give a more explicit expression of $\kappa_{\mc H}(\lambda)$, which is the content of the following theorem.
\begin{theorem}
\label{teo:2}
Let $\lambda$ be a simple eigenvalue of a real Hamiltonian matrix $Q$ with corresponding right and left eigenvectors $x$ and $y$ normalized to have
\begin{equation}
\label{p:1tris}
\|x\|_2=\|y\|_2 = 1\;, \qquad \Im(y^*Jx)\, \Re(y^*Jx) = \mathrm{Trace}(\Re(yx^*)^T \Im(yx^*)) \;.
\end{equation}
Given any real Hamiltonian matrix $E$ with $\|E\|_F = 1$, let $\lambda_E(t)$ be an eigenvalue of $Q+tE$ converging to $\lambda$ as $t\to 0$. Then,
\begin{equation}
\label{p:10}
\begin{split}
|\dot\lambda_E(0)| & \le \max\left\{\left|\frac{y^*Hx}{y^*x}\right| \colon  \|H\|_F = 1,\, H \in \mc H \,\right\} \\ & = \max\left\{\frac{\|\Re(yx^*)|_\mc H\|_F}{|y^*x|}; \frac{\|\Im(yx^*)|_\mc H\|_F\}}{|y^*x|}\right\}\;.\end{split}
\end{equation}
Moreover:
\\
1)  if $\|\Re(yx^*)|_\mc H\|_F>\|\Im(yx^*)|_\mc H\|_F$ then
\[
\dot\lambda_E(0) = \frac{\|\Re(yx^*)|_\mc H\|_F}{y^*x} \quad \mathit{for} \quad E=\pm\frac{\Re(yx^*)|_{\mc H}}{\|\Re(yx^*)|_\mc H)\|_F}\;;
\]
\\
2)  if $\|\Re(yx^*)|_\mc H\|_F<\|\Im(yx^*)|_\mc H\|_F$ then
\[
\dot\lambda_E(0) = \frac{\|\Im(yx^*)|_\mc H\|_F}{y^*x} \quad \mathit{for} \quad E=\pm\frac{\Im(yx^*)|_{\mc H}}{\|\Im(yx^*)|_\mc H)\|_F}\;;
\]
\\
3)  if $\|\Re(yx^*)|_\mc H\|_F=\|\Im(yx^*)|_\mc H\|_F$ then
\[
\dot\lambda_E(0) = \frac{\|\Re(yx^*)|_\mc H\|_F}{y^*x} \quad \mathit{for}\; \mathit{any} \quad E=\frac{\Re(yx^*)|_{\mc H}\cos\theta + \Im(yx^*)|_{\mc H} \sin\theta}{\|\Re(yx^*)|_\mc H\|_F}\;, 
\]
with $\theta\in [0,2\pi]$.
\end{theorem}
\proof Our target is to compute 
\[
\mathrm{argmax}\left\{|y^*H x| \colon \;H\in\mc H\;,\; \|H\|_F^2=1\right\}.
\]
As before, we introduce the real matrices $\xi = \Re(yx^*)$,  $\eta = \Im(yx^*)$. By arguing as in the proof of Theorem \ref{teo:1}, we find that the argmax have to be of the following form,
\begin{equation}
\label{abr}
H = \alpha (\xi+J\xi^TJ)  + \beta (\eta + J\eta^TJ)\;.
\end{equation}
As $\lan\xi+J\xi^TJ, \eta + J\eta^TJ \ran = 2\lan\xi,\eta\ran+2\lan\xi,J\ran\lan\eta,J\ran$ and recalling $y^*Jx =\lan \xi,J\ran -\rmi \lan\eta,J\ran$, by \eqref{p:1tris} we have,
\[
\|\xi\|_F^2+\|\eta\|_F^2 = 1\;, \qquad \lan \xi+J\xi^TJ, \eta + J\eta^TJ \ran = 0\;.
\]
Whence, a straightforward computation shows that if $H= \alpha (\xi+J\xi^TJ)  + \beta (\eta + J\eta^TJ)$ then
\[
| y^*H x|^2 = D_\xi^2  \alpha^2 + D_\eta^2\beta^2\;,\quad \|H\|_F^2= 2D_\xi\alpha^2 + 2D_\eta \beta^2\;,
\]
where
\[
\begin{split}
D_\xi  & := \frac 12 \|\xi+J\xi^TJ\|_F^2 = \|\xi\|_F^2+\lan\xi,J\ran^2\;, \\ D_\eta & := \frac 12 \| \eta + J\eta^TJ\|_F^2 = \|\eta\|_F^2+\lan\eta,J\ran^2\;.
\end{split}
\]
Moreover, as $\Re(yx^*)|_\mc H$ [resp.\ $\Im(yx^*|_\mc H)$] is given by \eqref{abr} for $(\alpha,\beta) = (1/2,0)$ [resp.\ $(\alpha,\beta)=(0,1/2)$], it follows that
\[
\|\Re(yx^*)|_\mc H\|_F^2 = \frac{D_\xi}2 \;, \quad \|\Im(yx^*)|_\mc H\|_F^2 = \frac{D_\eta}2\;.
\]
If both $D_\xi$ and $D_\eta$ are positive the (constrained) maximum of $| y^*H x|^2$ is reached for $(\alpha^2,\beta^2) = ((2D_\xi)^{-1},0)$ if $D_\xi>D_\eta$ and for $(\alpha^2,\beta^2)=(0,(2D_\eta)^{-1})$  if $D_\xi<D_\eta$. In the critical case $D_\xi=D_\eta=:D$, the function $| y^*Hx|^2$ is constant and any pair $(\alpha^2,\beta^2)$ satisfying the constraint  $2D(\alpha^2 +\beta^2)=1$ is allowed. Therefore, the maximal perturbation is obtained for 
\[
\begin{cases} {\displaystyle E_\xi =\pm\frac{ \xi+J\xi^TJ}{\sqrt{2D_\xi}} }& \text{if } D_\xi > D_\eta\;, \\ \\ {\displaystyle E_\eta = \pm\frac{ \eta+J\eta^TJ}{\sqrt{2D_\eta}} } & \text{if } D_\xi < D_\eta\;, \end{cases}
\]
while, in the critical case $D_\xi=D_\eta=:D$, for
\[
E_\theta = \frac{ (\xi+J\xi^TJ)\cos\theta  + (\eta + J\eta^TJ)\sin\theta}{\sqrt{2D}}\;, \quad\theta\in [0,\pi]\;.
\]
On the other hand, if $D_\eta=0$ [resp.\ $D_\xi=0$] then $\alpha^2=(2D_\xi)^{-1}$ [resp.\ $\beta^2 = (2D_\eta)^{-1}$] and $H=E_\xi$ [resp.\ $H=E_\eta$] holds trivially. In all cases, 
\[
| y^*E_\xi x| = \sqrt \frac{D_\xi}2\;,\quad | y^*E_\eta x| = \sqrt\frac{D_\eta}2\;, \quad | y^*E_\theta x| = \sqrt\frac D2\;.
\]
The theorem is thus proved. 
\qed
\begin{remark}
\label{rem:p3}
We observe that $yx^*$ is Hamiltonian if and only if $\Im(yx^*)|_\mc H=0$ and $\Re(yx^*) = \Re(yx^*)|_\mc H$. Likewise, $yx^*$ is skew-Hamiltonian if and only if $\Re(yx^*)|_\mc H=0$ and $\Im(yx^*) = \Im(yx^*)|_\mc H$. Therefore, by Theorem \ref{teo:2}, 
\[
\kappa_\mc H(\lambda) = \begin{cases} {\displaystyle \frac{\|\Re(yx^*)\|_F}{|y^*x|}}  & \mathrm{if}\;\;yx^*\;\;\mbox{is Hamiltonian}\;, \\ {\displaystyle  \frac{\|\Im(yx^*)\|_F}{|y^*x|} }& \mathrm{if}\;\;yx^*\;\;\mbox{is skew-Hamiltonian}\;. \end{cases}
\]
By \cite[Theorem 2.1]{BK03} it follows that $\kappa_\mc H(\lambda)$ is equal to the real (unstructured) condition number. Moreover, by Proposition \ref{prop:p1} this occurs only if  $\lambda$ is a purely imaginary eigenvalue.
\end{remark}
\begin{remark}
\label{rem:p4}
If is worthwhile to remark that if the eigenvalue $\lambda$ is real, since the right and left eigenvectors are also real, then \eqref{p:1} and \eqref{p:1tris} are trivially satisfied and $\kappa_{\mc H}(\lambda) = 
\kappa_{\mc H_{\bb C}}(\lambda)$.
\end{remark}

\section{Numerical tests}
\label{sec:4}

\noindent{\bf Example 4.1.} Consider the Hamiltonian matrix $H = \begin{pmatrix} H_1 & H_3 \\ H_2 & H_4 \end{pmatrix}\in \bb R^{4\times 4}$, where
\[
H_1=\begin{pmatrix} 0 & 2 \\ 0 & 0 \end{pmatrix}, \quad H_2=\begin{pmatrix} -0.1 & 0.1 \\ \phantom{2}0.1 & 0.1 \end{pmatrix}, \quad
H_3=\begin{pmatrix} 0 & \phantom{2}0.1 \\ 0.1 & -0.1 \end{pmatrix}.
\] 
It is well known that if $\lambda \in \bb C - (\bb R \cup \rm \rmi \bb R)$ is an eigenvalue of a real Hamiltonian matrix, then also $-\lambda$, $\conj\lambda$, $-\conj\lambda$ are eigenvalues, since the spectrum is symmetric with respect to the real and the imaginary axis.
According to Theorems \ref{teo:1} and \ref{teo:2}, one derives for the eigenvalues of $H$ the condition numbers reported in Table \ref{tab1}.
\begin{table}[tbh!]
\begin{center}
\begin{tabular}{c|c|c|c} 
$\lambda$ & $\kappa(\lambda)$ & $\kappa_{\mc H_{\bb C}}(\lambda)$ & $\kappa_\mc H(\lambda)$  \\
\hline
$\pm 0.3205 \pm \rmi 0.3126$ & $2.3513$ &$1.6631$ & $1.5851$ \\
\end{tabular}
\end{center}
\caption{Example 4.1. Unstructured, structured and real structured eigenvalue condition numbers.}\label{tab1}
\end{table}
\begin{figure}
\includegraphics[scale=0.60]{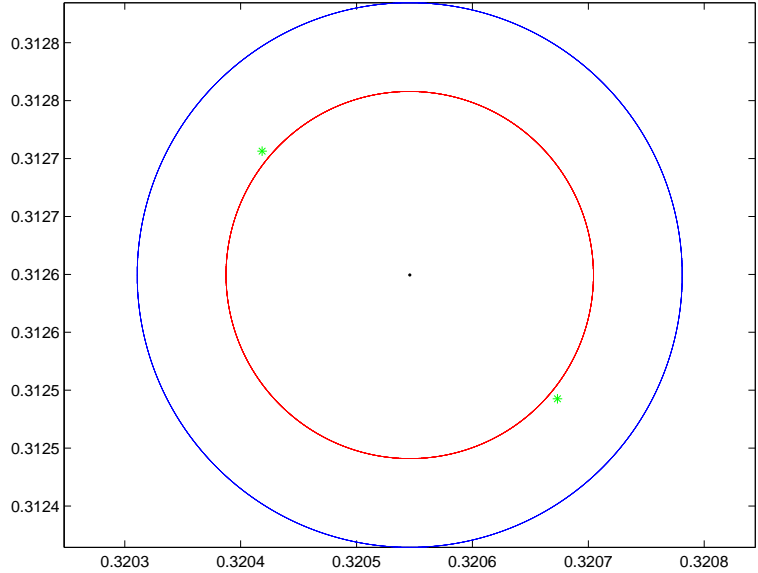}  
\caption{Example 4.1. Perturbations induced in the eigenvalue $\lambda$ by the complex structured maximal perturbations (green $\ast$),  the real structured maximal perturbations (red $-$) and the unstructured maximal perturbations (blu $-$).}\label{fig1}
\end{figure}
Let us focus on the eigenvalue $\lambda=0.3205 + \rmi 0.3126$, marked by a black point in Figure \ref{fig1}, and consider $\varepsilon$-norm perturbations, i.e. perturbation matrices of  Frobenius norm equal to $\varepsilon$, with  $\varepsilon = 10^{-4}$. 

The complex structured maximal perturbation matrices in (\ref{p:2}) shift  $\lambda$ to one of the points marked by a green star in Figure \ref{fig1}.  As a matter of fact, with respect to structured perturbations, the eigenvalue is not equally sensitive in any direction of the complex plane, since its eigenvectors have been normalized such that \eqref{p:1} is satisfied.

As for the real  structured  perturbations, holding  in this case the equality $\|\Re(yx^*)|_\mc H\|_F=\|\Im(yx^*)|_\mc H\|_F$, the maximal perturbation matrices are given by  case 3) in Theorem  \ref{teo:2} and they shift $\lambda$ to the boundary of the red circle of radius $\kappa_\mc H(\lambda) \varepsilon$. 

We point out that the unstructured worst-case effect perturbation matrices $\varepsilon \,\rme^{\rmi \theta}yx^*$, for  $\theta\in [0,2\pi]$, shift in their turn  $\lambda$ to the boundary of the blue circle of radius $\kappa(\lambda) \varepsilon$. Analogous results hold for the other three eigenvalues of $H$.

\medskip
\noindent{\bf Example 4.2.} Consider the same Hamiltonian test matrix as in Example 4.1, except for the block $H_3$ that is here multiplied by the factor $0.1$.
\begin{table}[tbh!]
\begin{center}
\begin{tabular}{c|c|c|c} 
$\lambda$ & $\kappa(\lambda)$ & $\kappa_{\mc H_{\bb C}}(\lambda)$ & $\kappa_\mc H(\lambda)$  \\
\hline
$\pm 0.1786 \pm \rmi 0.1771$ & $7.0263$ &$6.3769$ & $6.3184$ \\
\end{tabular}
\end{center}
\caption{Example 4.2. Unstructured, structured and real structured eigenvalue condition numbers.}\label{tab2}
\end{table}
\begin{figure}
\includegraphics[scale=0.60]{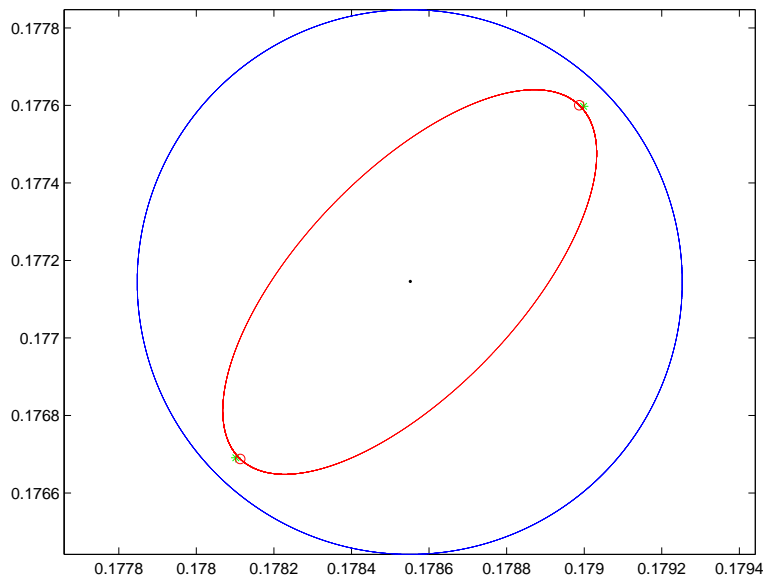}  
\caption{Example 4.2. Perturbations induced in $\lambda$ by the complex structured (green $\ast$), real structured (red $\circ$) and unstructured (blu $-$) maximal perturbations, and by the real structured  perturbations in (\ref{Eth}) (red $-$).}\label{fig2}
\end{figure}
We consider $\varepsilon$-norm perturbations in the matrix  and analyze the effects induced in $\lambda=0.1786  + \rmi 0.1771$.  The complex [real] structured maximal perturbation matrices shift  $\lambda$ to one of the points marked by a green star [red circlet] in Figure \ref{fig2}. 
Since $\|\Re(yx^*)|_\mc H\|_F<\|\Im(yx^*)|_\mc H\|_F$, case 2) in Theorem  \ref{teo:2} holds; nevertheless, for the sake of completeness,  let us take into account even the (non-maximal) real Hamiltonian matrices
\begin{equation}
\label{Eth}
 E_{\theta}=\frac{\Re(yx^*)|_{\mc H}\cos\theta + \Im(yx^*)|_{\mc H} \sin\theta}{\|\Re(yx^*)|_{\mc H}\cos\theta + \Im(yx^*)|_{\mc H} \sin\theta\|_F}\;,
\end{equation}
and notice that the perturbation matrices $\varepsilon \,E_{\theta}$, for  $\theta\in [0,2\pi]$,  have the effect of shifting $\lambda$ to the boundary of the  ellipse in Figure \ref{fig2}, whose major axis gives indeed the direction induced by  the worst-case effect perturbations.

\section{Conclusions}
\label{sec:5}
We have derived simple expressions for structured eigenvalue condition number and for the  worst-case effect perturbation in case of complex and real Hamiltonian matrices, which turns out to be a rank-2 complex matrix in the former and a rank-4 real matrix in the latter. In Proposition \ref{prop:p1} and Remark \ref{rem:p3}, we showed that the structured and unstructured eigenvalue condition numbers are equal in case of purely imaginary eigenvalues,  the unstructured worst-case effect perturbations resulting to be Hamiltonian. On the other hand, numerical tests show that the structured perturbation analysis may be of some interest for eigenvalues not lying in the imaginary axis; see, e.g., the tests reported in Section \ref{sec:4}.


\begin{thebibliography}{99}

\bibitem [BK05]{BK05}
P.~Benner, D.~Kressner, New Hamiltonian Eigensolvers with Applications in Control, Proceedings of 44th IEEE Conference on Decision and European Control Conference ECC 2005, pp.~655--6556. 

\bibitem [BKM05]{BKM05}
P.~Benner, D.~Kressner, and V.~Mehrmann, Skew-Hamiltonian and Hamiltonian Eigenvalue Problems: Theory, Algorithms and Applications. Z. Drmac, M. Marusic, and Z. Tutek, editors, Proceedings of the Conference on Applied Mathematics and Scientific Computing 2003, Springer-Verlag, 2005, pp.~3--39. 

\bibitem [BMX98]{BMX98}
P.~Benner, V.~Mehrmann, and H.~Xu, A numerically stable, structure preserving method for computing the eigenvalues of real Hamiltonian or symplectic pencils, Numer. Math., 78 (1998), pp.~329--358.

\bibitem[B87]{B87}
J.~R.~Bunch, The weak and strong stability of algorithms in numerical linear algebra. Linear Algebra Appl., 88/89 (1987), pp.~49Ð66.

\bibitem [BK03] {BK03}
R.~Byers, D.~Kressner, On the condition of a complex eigenvalue under real perturbations,  BIT, 43 (2003), pp.~1Ð18. 

\bibitem[KKT06] {KKT}
M.~Karow, D.~Kressner, and F.~Tisseur, Structured eigenvalue condition numbers,
SIAM J. Matrix Anal. Appl., 28 (2006), pp.~1052--1068.

\bibitem[NP06]{NP06}
S.~Noschese, L.~Pasquini, Eigenvalue condition cumbers: zero-structured versus traditional, J. Comput. Appl. Math., 185 (2006), pp.~174--189.

\bibitem[NP07]{NP07}  
S.~Noschese, L.~Pasquini, Eigenvalue patterned condition 
numbers: Toeplitz and Hankel cases,  J. Comput. Appl. Math., 206 (2007), pp.~615--624. 

\bibitem[R06] {R06}
S.~M.~Rump,  Eigenvalues, pseudospectrum and structured
perturbations, Linear Algebra Appl., 413 (2006), pp.~567--593.

\bibitem[T03]{T03} 
F.~Tisseur, A chart of backward errors for singly and doubly structured eigen-
value problems, SIAM J. Matrix Anal. Appl., 24 (2003), pp.~877--897.


\bibitem[W65]{W65} 
J.~H.~ Wilkinson, The algebraic eigenvalue problem, Oxford University Press, 1965.


\end{thebibliography}
\end{document}